\documentclass[11pt]{article}
\makeindex
\usepackage{textcomp}
\usepackage{epsfig, color}
\usepackage{amssymb,amsmath,latexsym}
\usepackage{epstopdf}
\usepackage{float}

\usepackage{stmaryrd}
\newcommand \brk[1]
{
\llbracket #1\rrbracket
}

\def \R {{\mathbb R}}

\newcommand \spann[1]
{
	\langle #1 \rangle
}

\def \EX {{\cal EX}}

\def \setp {{\cal P}}

\usepackage{tikz}
\usetikzlibrary{backgrounds}
\usetikzlibrary{arrows}
\usetikzlibrary{shapes,shapes.geometric,shapes.misc}
\pgfdeclarelayer{edgelayer}
\pgfdeclarelayer{nodelayer}
\pgfsetlayers{background,edgelayer,nodelayer,main}
\tikzstyle{none}=[inner sep=0mm]
\tikzstyle{every loop}=[]

\tikzstyle{dotted}=[dash pattern=on \pgflinewidth off 2pt]
\tikzstyle{dashed}=[dash pattern=on 3pt off 3pt]

\newcommand \tikzp[2]
{
\begin{center}
\begin{tikzpicture}[scale=#1]
#2
\end{tikzpicture}
\end{center}
}

\tikzstyle{new style 0}=[fill=black, draw=black, shape=circle]
\tikzstyle{red style 1}=[fill=red, draw=black, shape=circle]
\tikzstyle{blue style 2}=[fill=blue, draw=black, shape=circle]
\tikzstyle{white style 4}=[fill=white, draw=black, shape=circle]
\tikzstyle{bklack style 5}=[fill=black, draw=black, shape=rectangle]
\tikzstyle{red style 3}=[fill=red, draw=black, shape=rectangle]
\tikzstyle{yellow style 7}=[fill=yellow, draw=black, shape=rectangle]
\tikzstyle{new style 8}=[fill={rgb,255: red,0; green,132; blue,0}, draw={rgb,255: red,0; green,131; blue,0}, shape=circle]

\tikzstyle{new edge style 0}=[-]
\tikzstyle{new edge style 1}=[-, draw=red]
\tikzstyle{new edge style 2}=[-, draw=blue]
\tikzstyle{new edge style 3}=[-, draw={rgb,255: red,0; green,156; blue,0}]

\tikzstyle{cblue}=[circle, draw, thin,fill=blue!20, scale=0.5]
\newtheorem{lem}{Lemma}[section]
\newtheorem{theo}{Theorem}[section]
\newtheorem{pro}{Proposition}[section]

\newtheorem{cor}{Corollary}[section]

\newtheorem{rem}{Remark}[section]

\newcommand\blue[1] {{\color{blue} #1}}

\newcommand \lemma[2]  %%%used here
{\begin{lem}
\label{#1} #2
\end{lem}
}

\newcommand \prop[2]
{\begin{pro}
\label{#1} #2
\end{pro}
}

\newcommand{\proof}{{\noindent {\em Proof}.\quad}\setcounter{countclaim}{0}\setcounter{countcase}{0}}

\newcommand{\proofend}{{\hfill$\Box$}}

\newcounter{countcase}

\newcounter{countclaim}
\def\inclaim{\addtocounter{countclaim}{1}
{\noindent {\bf Claim \thecountclaim}: }}

\newcounter{countfig}

\newcommand{\beeq}{\begin{equation}}
\newcommand{\eneq}{\end{equation}}

\newcommand{\beeqn}{\begin{eqnarray*}}
\newcommand{\eneqn}{\end{eqnarray*}}

\def \setf{{\cal F}}
\def \setg{{\cal G}}

\def \setp{{\cal P}}
\def \sets{{\cal S}}

\def \sett{{\cal T}}

\def \ST {{\cal ST}}
\def \CT {{\cal CG}}

\def \tutte {{\bf T}}

\def \iff {if and only if }

\begin{document}

\baselineskip 0.6 cm

\title
{Proving identities on weight polynomials of tiered trees via Tutte polynomials}

\author{
Fengming Dong\thanks{fengming.dong@nie.edu.sg
and donggraph@163.com.}
\\
\small National Institute of Education,
Nanyang Technological University, Singapore
\\
Sherry H.F. Yan\thanks{Corresponding author.
Email: hfy@zjnu.cn.}
\\
\small Department of Mathematics,
Zhejiang Normal University, P.R. China
}

\date{}

\maketitle{}

\begin{abstract}

 A {\it tiered graph} $G=(V,E)$ with $m $ tiers is a simple graph
with $V\subseteq \brk{n}$, where $\brk{n}=\{1,2,\cdots,n\}$, and with
a surjective map $t$ from $V$ to $\brk{m}$
such that if $v$ is a vertex adjacent to $v'$ in $G$
with $v>v'$, then $t(v) >t(v')$.
 For any ordered partition $p=(p_1,p_2,\cdots,p_m)$ of $n$,
  let $\sett_p$ denote  the set of tiered trees
with vertex set $\brk{n}$ and
with a map $t: \brk{n}\rightarrow \brk{m}$
such that $|t^{-1}(i)|=p_i$   for all $i=1,2,\ldots,m$.
 For any $T\in \sett_p$,
let $K_T$ denote the complete
tiered graph   whose vertex set    and  tiering map  are  the same as  those of  $T$.
If the edges of $K_T$ are ordered lexicographically by their endpoints,
then the weight $w(T)$ of $T$ is the external activity of $T$ in $K_T$,
i.e., the number of edges
$e\in E(K_{T})\setminus E(T)$  such that $e$ is the
least element in the unique cycle determined by $T\cup e$.
\label{r2-1}
Let $P_p(q)=\sum_{T\in \sett_{p}}q^{w(T)}$.
Dugan, Glennon, Gunnells and Steingr\'imsson
[J. Combin. Theory, Ser. A 164 (2019) pp. 24-49]
asked for an elementary proof of the identity
$P_p(q)=P_{\pi(p)}(q)$
for any permutation $\pi$ of $1,2,\cdots,m$,
where $\pi(p)=(p_{\pi(1)},p_{\pi(2)},\cdots,p_{\pi(m)})$.
 In this article,
	we will prove an extension of this identity by applying Tutte polynomials.
Furthermore, we also provide a proof
of the  identity $P_{(1,p_1,p_2)}(q)=P_{(p_1+1,p_2+1)}(q)$
via Tutte polynomials.
\end{abstract}



%\tableofcontents

\section{Introduction}

\subsection{Tiered graphs}

A {\it tiered graph} $G=(V,E)$ with $m \ge 2$ tiers is a simple graph
with $V\subseteq \brk{n}$, where $\brk{n}=\{1,2,\cdots,n\}$, and with
a surjective map $t$ from $V$ to $\brk{m}$
such that if $v$ is a vertex adjacent to $v'$ in $G$
with $v>v'$, then $t(v) >t(v')$.
We call $t$ a {\it tiering map} of $G$ and
say $G$ is {\it tiered} by $t$.
If a tiered graph is a tree, it
is called a {\it tiered tree}.
The concept of tiered trees was introduced by %Dugan et al.
Dugan, Glennon, Gunnells and Steingr\'imsson
\cite{dugan2019tiered} who generalized
the concept of intransitive trees
(also called alternating trees)
introduced by Postnikov \cite{postnikov1997intransitive},  the latter of which have exactly
two tiers. \label{r2-2}

 As  mentioned in \cite{dugan2019tiered},
tiered trees naturally arise in two unrelated
geometric counting problems \cite{godsil2013algebraic}:
counting absolutely irreducible representations of the supernova quivers
and counting certain torus orbits on partial flag varieties of type A
over finite fields, namely those orbits with trivial stabilizers.

If $G$ is a tiered graph with $m$ tiers and a  tiering map $t$, let  $V_i(G)=\{v\in V(G): t(v)=i\}$ for all $i\in \brk{m}$.
Clearly, each set $V_i(G)$ is an independent set in $G$,
and for any $v\in V_i(G)$ and $v'\in V_j(G)$,
if $vv'\in E(G)$, then $v<v'$ \iff $i<j$.

\iffalse
For any ordered partition $p=(p_1,p_2,\cdots,p_m)$ of
$n$ (i.e., $p_i\ge 1$ and $p_1+p_2+\cdots+p_m=n$),
let $\setg_p$ (resp. $\sett_p$) be the set of tiered graphs
(resp. tiered trees)
with vertex set $\brk{n}$
and $|V_i(G)|=p_i$ for $i=1,2,\cdots,m$.
Obviously, for any graph $G$ with $V(G)=\brk{n}$,
$G\in \setg_p$ \iff $V(G)$ has a partition
$V_1, V_2,\cdots, V_m$ such that the following properties hold:

(a). for $i=1,2,\cdots,m$,
$V_i$ is an independent set of $G$ with $|V_i|=p_i$; and

(b). for any edge $vv'\in E(G)$, where $v\in V_i$ and $v'\in V_j$,
$v<v'$ \iff $i<j$.

There may exist graphs $G_1,G_2\in \setg_p$
with $V(G_1)=V(G_2)$ and $E(G_1)=E(G_2)$, but
$V_i(G_1)\ne V_i(G_2)$ for some $i\in \brk{m}$.
\fi
In this article, we assume that
two tiered graphs  $G_1$ and $G_2$  are {\it different}
if either $G_1$ and $G_2$ have different edge sets or  $V_i(G_1)\ne V_i(G_2)$ for some $i$.

When we draw a tiered graph $G$, vertices in the same
set $V_i(G)$ are put at the same height
and vertices in $V_j(G)$ are at a higher level
than vertices in $V_i(G)$ whenever $i<j$.
For example, $\sett_{(2,2)}$ contains five tiered trees
(see \cite{dugan2019tiered})
as shown in Figure~\ref{fig1},
where $\sett_{(2,2)}$ is the set of tiered trees with vertex set $\brk{4}$
and two tier sets of size $2$ each. %$|V_1(T)|=|V_2(T_1)|=2$.

\begin{figure}[H]
\centering

\tikzp{1}
{
%%%%%%%%%%%%%%%%%%% figure (a)
\foreach \place/\z in
{{(4,0)/1}, {(2,0)/2},{(1,0.5)/3},{(3,0.5)/4},
{(4,1.5)/5},{(2,1.5)/6}, {(3,2)/7},{(1,2)/8}}
\node[cblue] (w\z) at \place {};

\filldraw[black] (w1) circle (0pt) node[anchor=south]
{\small $1$};
\filldraw[black] (w2) circle (0pt) node[anchor=south]
{\small $2$};
\filldraw[black] (w3) circle (0pt) node[anchor=south]
{\small $3$};
\filldraw[black] (w4) circle (0pt) node[anchor=south]
{\small $4$};
\filldraw[black] (w5) circle (0pt) node[anchor=south]
{\small $1$};
\filldraw[black] (w6) circle (0pt) node[anchor=south]
{\small $2$};
\filldraw[black] (w7) circle (0pt) node[anchor=south]
{\small $3$};
\filldraw[black] (w8) circle (0pt) node[anchor=south]
{\small $4$};

\draw[black] (w3) -- (w2) -- (w4) -- (w1);

\draw[black] (w8) -- (w6) -- (w7) -- (w5);

\filldraw[black] (2.5,-0.25) circle (0pt) node[anchor=north] {(a)};

%%%%%%%%%%%%%%%%%%%%%%%%% figure (b)
\foreach \place/\z in
{{(7,0)/1}, {(9,0)/2},{(6,0.5)/3},{(8,0.5)/4},
{(7,1.5)/5},{(9,1.5)/6}, {(8,2)/7},{(6,2)/8}}
\node[cblue] (u\z) at \place {};

\filldraw[black] (u1) circle (0pt) node[anchor=south]
{\small $1$};
\filldraw[black] (u2) circle (0pt) node[anchor=south]
{\small $2$};
\filldraw[black] (u3) circle (0pt) node[anchor=south]
{\small $3$};
\filldraw[black] (u4) circle (0pt) node[anchor=south]
{\small $4$};
\filldraw[black] (u5) circle (0pt) node[anchor=south]
{\small $1$};
\filldraw[black] (u6) circle (0pt) node[anchor=south]
{\small $2$};
\filldraw[black] (u7) circle (0pt) node[anchor=south]
{\small $3$};
\filldraw[black] (u8) circle (0pt) node[anchor=south]
{\small $4$};

\draw[black] (u3) -- (u1) -- (u4) -- (u2);

\draw[black] (u8) -- (u5) -- (u7) -- (u6);

\filldraw[black] (7.5,-0.25) circle (0pt) node[anchor=north] {(b)};

%%%%%%%%%%%%%figure (c) below
\foreach \place/\z in
{{(12,0)/1}, {(11,0.5)/2},{(14,0)/3},{(13,0.5)/4}}
\node[cblue] (v\z) at \place {};

\filldraw[black] (v1) circle (0pt) node[anchor=south]
{\small $1$};
\filldraw[black] (v2) circle (0pt) node[anchor=south]
{\small $2$};
\filldraw[black] (v3) circle (0pt) node[anchor=south]
{\small $3$};
\filldraw[black] (v4) circle (0pt) node[anchor=south]
{\small $4$};

\draw[black] (v2) -- (v1) -- (v4) -- (v3);

\filldraw[black] (12.5,-0.25) circle (0pt) node[anchor=north] {(c)};
}

\caption{Tiered trees in $\sett_{(2,2)}$}
\label{fig1}
\end{figure}\label{r1-1}

A non-tiered graph
may be not simple and may have parallel edges and loops,
	unless  otherwise stated.
	For any graph $G$,
	let $V(G)$ and $E(G)$ denote the vertex set and edge set of $G$.
	For any  $E_0\subseteq E(G)$,
	let $G\spann{E_0}$ (or $(V(G), E_0)$)
	denote the spanning subgraph of $G$ with edge set $E_0$.
	Let $\ST(G)$ denote the set of spanning trees of $G$.
	A component of $G$ is referred to  a connected  component of $G$.
	For any graphs $G$ and $H$,
		let $G\cup H$ denote the graph with vertex set $V(G)\cup V(H)$ and edge set $E(G)\uplus E(H)$,
	 where the symbol $\uplus$ indicates that each pair of edges
	$e_1\in E(G)$ and $e_2\in E(H)$ are treated as
	 	different edges in $G\cup H$.

	For any non-empty finite set $U$ of positive integers
	and an ordered partition $p=(p_1,p_2,\cdots,p_m)$ of $|U|$, where $m\ge 2$,
	let $\setg_{U,p}$ denote the set of tiered graphs $G$ with  $V(G)=U$,   $m$ tiers and
	$|V_i(G)|=p_i$ for all $i=1,2,\cdots,m$,
	and let $\sett_{U,p}$ (resp. $\setf_{U,p}$) be the set of trees
	(resp. forests) in $\setg_{U,p}$.
	
	A tiered graph $G\in \setg_{U,p}$
	is said to be {\it complete}
	if for any $1\le i<j\le m$,
	$u\in V_i(G)$ and $v\in V_j(G)$,
	$uv\in E(G)$ \iff $u<v$.
	Let $\CT_{U,p}$ denote the set of compete tiered graphs in $\setg_{U,p}$
	and let $\CT^c_{U,p}$ denote the set of connected graphs $G\in \CT_{U,p}$.
	In the case that $U=\brk{n}$ for some positive integer $n$, $\setg_{U,p}$, $\sett_{U,p}$,
	$\CT_{U,p}$ and $\CT^c_{U,p}$ are simply written as
	$\setg_{p}$, $\sett_{p}$, $\CT_{p}$ and $\CT^c_{p}$ respectively.

	 For a graph $H$ (possibly not tiered) with  $U\subseteq V(H)$,
let $\setg_{U,p}(H)$ (resp. $\CT_{U,p}(H)$)
denote the set of graphs $H\cup Q$ for $Q\in \setg_{U,p}$ (resp. $Q\in \CT_{U,p}$).
Let $\CT^c_{U,p}(H)$ denote   the set of connected graphs in $\CT_{U,p}(H)$,
and let $\sett_{U,p}(H)$  denote the set of trees in $\setg_{U,p}(H)$.
Obviously, $\sett_{U,p}(H)=\emptyset$ if $H$ is not acyclic.

\subsection{Two identities  on   weight polynomials  of tiered trees
}

%Let $\sett=\bigcup_{p} \sett_p$.
 For any tiered tree $T$,  let $w(T)$ denote the weight of $T$
defined below \cite{dugan2019tiered}:
\begin{enumerate}
\item  if $|V(T)| =1$, then $w(T) =0$; and

\item  if $|V(T)| >1$, let $v=\min\{u:u\in V(T)\}$.    Suppose that
$T_1,\cdots, T_\ell$ are  the   components of the forest obtained by deleting $v$ from $T$ and
$u_i$ is the vertex in $T_i$
which is adjacent to $v$ for $i=1,2,\cdots,\ell$.
Let
$$
w(T) =\sum_{i=1}^\ell (w_i + w(T_i)),
$$
where $w_i$ is the cardinality of the set
$\{u_j\in V(T_i): u_j<u_i, t(u_j)>t(v)\}.
$
\end{enumerate}

For example, it can be checked easily that
all tiered trees in Figure~\ref{fig1} have weight $0$,
except the bottom one in Figure~\ref{fig1} (a),
which has weight $1$.
More examples are provided in \cite{dugan2019tiered}.

For any ordered partition $p=(p_1,p_2,\cdots,p_m)$ of $n$,
the weight polynomial for tiered trees in $\sett_p$
is defined
as follows:
\begin{equation}\label{eq1}
P_p(q)=\sum_{T\in \sett_p}q^{w(T)}.
\end{equation}

Dugan, Glennon, Gunnells and Steingr\'imsson \cite{dugan2019tiered}
remarked that, by applying the geometric results in \cite{gunnells2018torus},
one can prove that
the polynomials $P_p(q)$
depend only on the partition determined by the tier type $p$,
and not the order of its parts
(i.e., Theorem~\ref{main-th}).
They asked for an elementary proof of this result.
Yan et al. \cite{yan2020bijective} recently provided
a combinatorial proof of Theorem~\ref{main-th}
for the case $q=1$.

\begin{theo}[\cite{dugan2019tiered}]\label{main-th}
For any ordered partition $p=(p_1,p_2,\cdots,p_m)$
of $n$ and any permutation $\pi$ of $1,2,\cdots,m$,
 we have $P_p(q)=P_{\pi(p)}(q)$,
where $\pi(p)=(p_{\pi(1)}, p_{\pi(2)},\cdots,p_{\pi(m)})$.
\end{theo}

Another identity on $P_p(q)$ is
$P_{(1,p_1,p_2)}(q)=P_{(p_1+1,p_2+1)}(q)$
(i.e., the result of Theorem~\ref{main-th4}).
For example,
$P_{(1,1,1)}(q)=P_{(2,2)}(q)=q+4$,
$P_{(1,1,2)}(q)=P_{(2,3)}(q)=q^2 + 5q + 11$,
and $P_{(1,2,2)}(q)=P_{(3,3)}(q)=q^4 + 6q^3 + 22q^2 + 51q + 66$.
More examples for this identity are provided in \cite{dugan2019tiered}.
In general, this identity is explained in \cite{dugan2019tiered} by
counting the number of maximal $T$-orbits in $G(k, n)$,
where $G(k, n)$ is the Grassmannian of $k$ planes in
${\mathbb F}^n_q$.

\begin{theo}[\cite{dugan2019tiered}]
\label{main-th4}
For any positive integers $p_1$ and $p_2$,
$P_{(1,p_1,p_2)}(q)=P_{(p_1+1,p_2+1)}(q)$ holds.
\end{theo}

In this article, we will provide elementary proofs for Theorems~\ref{main-th} and \ref{main-th4}
via Tutte polynomials.

\subsection{An extension of Theorem~\ref{main-th}
	\label{sec1-3}}

For any graph $G=(V,E)$,
let $\tutte_G(x,y)$ denote the Tutte polynomial of $G$
defined recursively by the operations of
deletion and contraction \cite{tutte2009ring}:
\begin{equation}\label{eq2}
\tutte_G(x,y)=
\left \{
\begin{array}{ll}
1, &\mbox{if }E=\emptyset;\\
y\tutte_{G\backslash e}(x,y), &\mbox{if } e \mbox{ is a loop of }G;\\
x\tutte_{G\slash e}(x,y), &\mbox{if } e \mbox{ is a bridge of }G;\\
\tutte_{G\slash e}(x,y)+\tutte_{G\backslash e}(x,y),
&\mbox{if }e \mbox{ is neither a bridge nor a loop of } G,
\end{array}
\right.
\end{equation}
where $G\backslash e$ (resp. $G/e$)
is the graph obtained from $G$ by removing (resp. contracting) $e$.

For any ordered partition $p=(p_1,p_2,\cdots,p_m)$ of $n$,
recall that $\CT_p$  denotes  the set of complete tiered graphs  in $\setg_p$, and
$\CT^c_p$ is the set of connected graphs in $\CT_p$.
Dugan, Glennon, Gunnells and Steingr\'imsson \cite[Theorem 2.9]{dugan2019tiered} actually established   a connection
between $P_p(q)$ and
Tutte polynomials:
\begin{equation}\label{eq1-4}
P_p(q) = \sum_{G\in \CT^c_p} \tutte_G(1,q).
\end{equation}
 By (\ref{eq1-4}), in order to prove Theorem~\ref{main-th},
it suffices to show that for each $r\in \brk{m-1}$,
\begin{equation}\label{eq1-4-1}
\sum\limits_{G\in \CT^c_p} \tutte_G(1,q)
=\sum\limits_{G\in \CT^c_{\pi_r(p)}}\tutte_G(1,q),
\end{equation}
where $\pi_r(p)$ is the ordered partition
of $n$ obtained from $p$ by exchanging
$p_r$ and $p_{r+1}$, i.e.,
 $\pi_r(p)=(p_1,\cdots,p_{r-1},p_{r+1}, p_r, p_{r+2},\cdots,p_m)$.
Let $r$ be a fixed number in $\brk{m-1}$.
Each graph $G\in \CT^c_p$ can be decomposed
into two subgraphs, temporarily denoted by $G_r$ and $G_r'$,  where
$G_r$ is the subgraph of $G$ induced
by $V_r(G)\cup V_{r+1}(G)$,
and  $G_r'$ is  the spanning subgraph  $G\spann{E(G)\setminus E(G_r)}$.
Let $\left(\CT^c_p\right )_r'$ be the set
of graphs $G'_r$ for all $G\in \CT^c_p$.
 For any $H\in \left(\CT^c_p\right )_r'$,
let
$\CT^c_{p}[H]$ denote the set of graphs 
$G\in \CT^c_p$  such that $G'_r$ is exactly
$H$.
Then,  (\ref{eq1-4-1}) holds if we can prove that for each $H\in \left(\CT^c_p\right )_r'$,
\begin{equation}\label{eq1-4-2}
	\sum_{G\in \CT^c_{p}[H]}
	\tutte_G(1,q)
	=\sum_{G\in \CT^c_{\pi_r(p)}[H]} \tutte_G(1,q).
	%\qquad \forall H\in \left(\CT^c_p\right )_r'.
\end{equation}
Note that $\CT^c_p[H]$  is actually the set of connected graphs $H\cup Q$, \label{r2-3}
where $Q$ is a complete tiered graph with
vertex set $V_r(H)\cup V_{r+1}(H)$
and two tiers corresponding to the ordered partition $(p_r,p_{r+1})$.
Motivated by this observation,
we will prove a result (i.e. Theorem~\ref{main-th3})
which is an extension of (\ref{eq1-4-2})
as graph $H$ in this theorem can be  any graph
  with $V(Q)\subseteq V(H)$.

\begin{theo}\label{main-th3}
Let $U$ be a non-empty finite set of positive
integers, let $p=(p_1,p_2)$ be an ordered partition of $|U|$
 and let  $p'=(p_2, p_1)$.
For any graph $H$ with $U\subseteq V(H)$,
the following identity holds:
\begin{equation}\label{eq1-6}
\sum_{G\in \CT^c_{U,p}(H)} \tutte_G(1,y)
=\sum_{G\in \CT^c_{U,p'}(H)} \tutte_G(1,y).
\end{equation}
\end{theo}

Observe that (\ref{eq1-4-2}) is the special case of Theorems~\ref{main-th3} when $H$ is a tiered graph and $U$ is an independent set of $H$.
By the above explanation on (\ref{eq1-4-1}) and (\ref{eq1-4-2}) and their relation
with Theorem~\ref{main-th},
 we have the following conclusion.

\begin{cor}\label{main-cor}
Theorem~\ref{main-th} follows from Theorem~\ref{main-th3}.
\end{cor}

The other sections in this article are arranged as follows:
\begin{itemize}
\label{r2-4}
\item
In Section~\ref{sec3},
we define the dual graph $G'$ in $\setg_{U,p'}$ of each graph $G$ in $\setg_{U,p}$,
and show that $F\rightarrow F'$ is a bijection from forests    $F\in \setf_{U,p}$ to forests  $F'\in \setf_{U,p'}$,
where $U$ is a non-empty finite set of
positive integers,   $p=(p_1, p_2)$ is an ordered partition of $|U|$, and $p'=(p_2, p_1)$.
 Relying on  the construction of  dual graphs, we  shall establish  a bijection $T\rightarrow T^*$
from $\bigcup_{G\in \CT_{U,p}(H)}\ST(G)$
to $\bigcup_{G\in \CT_{U,p'}(H)}\ST(G)$,
 where $H$ is a graph with $U\subseteq V(H)$.
This bijection will be applied in the next section to obtain some preliminary results for Tutte polynomials.

\item
In Section~\ref{sec5},   \label{r2-5}
we will establish some preliminary results
on external activities of spanning trees
in graphs $G\in \CT_{U,p}(H)$, where  $U\subseteq V(H)$   and $p=(p_1, p_2)$ is an ordered partition of $|U|$.

\item
Theorem~\ref{main-th3}  is then proved in Section~\ref{sec6},
while Theorem~\ref{main-th4} is verified
in Section~\ref{sec7}
by establishing an identity
on Tutte polynomials of graphs in
$\CT_{(p_1,1,p_2)}$ and $\CT_{(p_1+1,p_2+1)}$
(i.e., Lemma~\ref{le7-1}).
\end{itemize}

\section{Dual graphs
\label{sec3}}

\subsection{The dual graph of a tiered  graph with two tiers
	\label{ns2-1}}

For any connected graph  $G$
with $V(G)=\{x_1,x_2,\cdots,x_s\}$ and
a tiering map $t: V(G)\rightarrow \brk{2}$,
where $1\le x_1<x_2<\cdots <x_s$,
let $G'$ denote the tiered graph
with the tiering map $t':V(G')\rightarrow \brk{2}$ such that
$V(G')=V(G)$,
$t'(x_{r})=3-t(x_{s+1-r})$ for $r=1,2,\cdots,s$, and
$x_ix_j\in E(G')$
\iff $x_{s+1-i}x_{s+1-j}\in E(G)$ for all $1\le i<j\le s$.
$G'$ is called the {\it dual graph} of $G$.
\iffalse
\blue{By definition, $G'$ is obtained from $G$ by exchanging the two tiers of $G$
	and changing the labels of all vertices from $x_i$ to $x_{s+1-i}$.
When drawing a diagram for $G'$,
we reorder the vertices
so that the vertex labels in each tier  are increasing
from left to right.
}\fi
For example,  a tree $T$ and its dual graph $T'$ are shown in Figure~\ref{f8} (a) and (b), respectively.

\begin{figure}[h!]
\centering

\tikzp{0.95}
{
%%%%%%%%%%%%%%%%%%% figure (a)
\foreach \place/\z in
{{(0,0)/2},{(1,0)/3},{(2,0)/4},{(3,0)/6},
{(0.5,1)/5}, {(1.5,1)/7},{(2.5,1)/8}}
\node[cblue] (w\z) at \place {};

\filldraw[black] (w2) circle (0pt) node[anchor=north]
{\small $2$};
\filldraw[black] (w3) circle (0pt) node[anchor=north]
{\small $3$};
\filldraw[black] (w4) circle (0pt) node[anchor=north]
{\small $4$};
\filldraw[black] (w5) circle (0pt) node[anchor=south]
{\small $5$};
\filldraw[black] (w6) circle (0pt) node[anchor=north]
{\small $6$};
\filldraw[black] (w7) circle (0pt) node[anchor=south]
{\small $7$};
\filldraw[black] (w8) circle (0pt) node[anchor=south]
{\small $8$};

\draw[black] (w2) -- (w5) -- (w4) -- (w8) -- (w6) -- (w7);

\draw[black] (w5) -- (w3);

%%%insert (a)
%\filldraw[black] (2.5,-0.5) circle (0pt) node[anchor=north] {(a)};

%%%%%%%%%%%%%%%%%%%%%%%%% figure (b)
\foreach \place/\z in
{{(5,0)/2},{(6,0)/3},{(7,0)/5},
{(4.5,1)/4}, {(5.5,1)/6}, {(6.5,1)/7},{(7.5,1)/8}}
\node[cblue] (w\z) at \place {};

\filldraw[black] (w2) circle (0pt) node[anchor=north]
{\small $2$};
\filldraw[black] (w3) circle (0pt) node[anchor=north]
{\small $3$};
\filldraw[black] (w5) circle (0pt) node[anchor=north]
{\small $5$};
\filldraw[black] (w4) circle (0pt) node[anchor=south]
{\small $4$};
\filldraw[black] (w6) circle (0pt) node[anchor=south]
{\small $6$};
\filldraw[black] (w7) circle (0pt) node[anchor=south]
{\small $7$};
\filldraw[black] (w8) circle (0pt) node[anchor=south]
{\small $8$};

\draw[black] (w3) -- (w4) -- (w2) -- (w6) -- (w5) -- (w7);

\draw[black] (w5) -- (w8);

}

{}\hfill \hspace{0.3 cm}
 (a) \hspace{3 cm} (b) \hfill {}

\caption{$T$ and its dual graph $T'$
}
\label{f8}
\end{figure}

Note that the dual tiered graph $T'$ of
$T$ in Figure~\ref{f8} is also a tree, which is even isomorphic to $T$.
This property actually holds for all tiered trees $T$ with two tiers.

\lemma
{le2-1}
{For any connected tiered  graph $G$
	with two tiers,  the following properties hold:
\begin{enumerate}
\item $G'$ is isomorphic to $G$;
\item $|V_i(G)|=|V_{3-i}(G')|$ holds for $i=1,2$;
\item $G$ is the dual tiered graph of $G'$ and $G\rightarrow G'$ is an involution.
\end{enumerate}
}

\proof
 (ii)  and (iii) follows  directly from the construction of the dual graph $G'$ of $G$. Assume that $V(G)=\{x_1,x_2,\cdots,x_s\}$,
where $1\le x_1<x_2<\cdots <x_s$.  By definition, $V(G')=V(G)$.
 By the definition of $G'$,
the map $\phi:V(G')\rightarrow V(G)$ defined by $\phi(x_i)=x_{s+1-i}$ for all $i=1,2,\cdots,s$ is an isomorphism from $G$ to $G'$.
 Thus, (i) holds, completing the proof.
\proofend

By Lemma~\ref{le2-1} (i) and (ii), a diagram of $G'$ can be obtained
by flipping a diagram of $G$ vertically
and changing the labels of all vertices from $x_i$ to $x_{s+1-i}$. \label{r2-8}

Note that, if $G$ is the graph with one vertex $x_1$,
then $t^{-1}(i)=\emptyset$ for some $i\in \brk{2}$.
By definition, if $t^{-1}(i)=\emptyset$, then
 $(t')^{-1}(3-i)=\emptyset$.

 Now we proceed to consider the case of disconnected graphs $G$.  \label{r2-9}
Assume that $G$ has   components $G_1,G_2,\cdots,G_c$.
For $1\le i\le c$, let
$t_i$ be the tiering map of $G_i$ defined by
$t_i(u)=t(u)$ for all $u\in V(G_i)$.
By Lemma~\ref{le2-1} (i), the dual graph $G_i'$ of $G_i$ is connected
for all $i=1,2,\cdots,c$.
The {\it dual graph} of $G$, denoted by $G'$,  \label{r2-10}
is the graph consisting of
components $G'_1,G'_2,\cdots, G'_c$
with the tiering map $t'$  defined by
$t'(u)=t'_{r_u}(u)$ for all $u\in V(G)$,
where $r_u$ is the unique number in $\brk{c}$
with $u\in V(G_{r_u})$.
 For example, a tiered forest $F$ and its dual graph $F'$ are illustrated in Figure \ref{f17} (b) and (c), respectively.
By Lemma \ref{le2-1} and
the definition of dual graphs, one can easily verify %that the dual graph of a tiered graph has
the following properties on dual graphs.

\lemma
{le2-2}
{Let $G$ be a tiered graph with a tiering map
$t:V(G)\rightarrow \brk{2}$.
Then, the following properties hold:
\begin{enumerate}
\item $G'\cong G$;
\item $|V_i(G')|=|V_{3-i}(G)|$ holds for $i=1,2$;
\item $G$ is the dual tiered graph of $G'$; and
\item for any $u,v\in V(G)$ with $u<v$,
$u$ and $v$ are in the same component of $G$
\iff $u$ and $v$ are in the same component of $G'$.
\end{enumerate}
}

	For any non-empty finite set $U$ of positive integers
	and ordered partition $p=(p_1,p_2)$ of $|U|$,
	recall that
	$\setg_{U,p}$ is the set of
	tiered graphs $G$ with $V(G)=U$
	and two tiers such that
	$|V_i(G)|=p_i$ for $i=1,2$.
By Lemma~\ref{le2-2} (ii),
$G'\in \setg_{U,p'}$ whenever $G\in \setg_{U,p}$, where $p'=(p_2,p_1)$.
\label{r2-11}

Recall that $\setf_{U,p}$ denotes
the set of forests in $\setg_{U,p}$.
Thus, $F'\in \setg_{U,p'}$ for any
$F\in \setg_{U,p}$ by Lemma~\ref{le2-2}.  \label{r2-12}

\prop
{pr2-1}
{
For any non-empty finite set $U$ of positive integers
and any  ordered partition $p=(p_1,p_2)$ of $|U|$,
$F\rightarrow F'$  is a bijection from
$\setf_{U,p}$ to $\setf_{U,p'}$.
}

 \proof
As $(p')'=p$,
by Lemma~\ref{le2-2},   \label{r2-14}
$\setf_{U,p}\ne \emptyset$ \iff $\setf_{U,p'}\ne \emptyset$.
Now assume that $\setf_{U,p}\ne \emptyset$.
Let $\psi$ denote the mapping $F\rightarrow F'$,
i.e., $\psi(F)=F'$ for each $F\in \setf_{U,p}$.
We first show that $\psi$ is surjective from
$\setf_{U,p}$ to $\setf_{U,p'}$.

Let $Q$ be any graph in $\setf_{U,p'}$.
By Lemma~\ref{le2-2}, $Q'\in \setf_{U,p}$.
So $\psi(Q')=(Q')'$.
By Lemma~\ref{le2-2} (iii), $(Q')'$ is $Q$ itself.
Thus, $\psi(Q')=Q$, and so
$\psi$ is surjective from
$\setf_{U,p}$ to $\setf_{U,p'}$.

 Now we proceed to show that $\psi$ is an injection from
$\setf_{U,p}$ to $\setf_{U,p'}$. For any graphs $Q_1, Q_2\in \setf_{U,p}$, assume that $\psi(Q_1)=Q'_1=Q'_2=\psi(Q_2)$. Again by Lemma~\ref{le2-2} (iii), we have $Q_1=\psi(Q'_1)=\psi(Q'_2)=Q_2$. This implies that $\psi$ is an injection from
$\setf_{U,p}$ to $\setf_{U,p'}$ as desired.
Hence the result holds.
 \proofend

\subsection{Bijections between quasi-tiered trees and
their dual trees
\label{ns2-2}
}

In this subsection,
let $U$ be a non-empty finite set of positive integers,
let $p=(p_1,p_2)$ be an ordered partition of $|U|$,
 $p'=(p_2,p_1)$ and
let $H$ be a graph
 with $U\subseteq V(H)$.

For any $Q\in \setg_{U,p}$,  \label{r2-15}
 $H\cup Q$ is  called a {\em quasi-tiered graph} with
 $V_i(H\cup Q)=V_i(Q)$ for $i=1,2$.
  Thus, $V_1(H\cup Q)\cup V_2(H\cup Q)=U$.
Note that $H\cup Q$ may be not a tiered graph,
 as $H$ may not  be a  tiered graph.
If  $H\cup Q$ is a tree, then it is called a {\em quasi-tiered tree}.  \label{r1-4}
%Note that $H\cup Q$ is a graph andan edge in $H$ may be parallel to an edge in $Q$.
We assume that two quasi-tiered graphs $G_1$ and $G_2$
are {\it different} if either
 $E(G_1)\ne E(G_2)$
or $V_i(G_1)\ne V_i(G_2)$ for some $i\in \brk{2}$.
Thus, for any two quasi-tiered graphs $G_1$ and $G_2$ with
$V_i(G_1)\ne V_i(G_2)$ for some $i\in \brk{2}$,
$\ST(G_1)\cap \ST(G_2)=\emptyset$.
For example,
$T_1$ and $T_2$ in Figure~\ref{f15} are two
different quasi-tiered trees since $V_1(T_1)\ne V_1(T_2)$.

For any $E_0\subseteq E(H)$ and
$Q\in \setg_{U,p}$,
let $E_0\cup Q$ denote the graph
$H\spann{E_0}\cup Q$.
Recall that $\setg_{U,p}(H\spann{E_0})$
is the set of graphs $E_0\cup Q$
(i.e., $H\spann{E_0}\cup Q$),
where $Q\in \setg_{U,p}$,
and $\sett_{U,p}(H\spann{E_0})$
is the set of trees in
$\setg_{U,p}(H\spann{E_0})$.
It is more convenient to
write  $\sett_{U,p}(E_0)$
for $\sett_{U,p}(H\spann{E_0})$.
 To be more specific, we have   \label{r2-19}
\begin{equation}\label{set-p-E0}
\sett_{U,p}(E_0)=\{E_0\cup F: F\in \setf_{U,p}
\mbox{ and }E_0\cup F \mbox{ is a tree} \}.
\end{equation}
 For any tree  $T=E_0\cup F\in \sett_{U,p}(E_0)$,   \label{r2-22}
the {\it dual graph} of $T$ with respect to $E_0$,
denoted by $T^*$,
is defined to be the quasi-tiered graph $E_0\cup F'\in \setg_{U,p'}(E_0)$,
where $F'\in \setf_{U, p'}$ is the dual forest of $F$.
An example of $T=E_0\cup F\in \sett_{U,p}(E_0)$ and $T^*=E_0\cup F'$ is shown in Figure~\ref{f17}, where $U=\brk{10}$,
	$p=(6,4)$,
$V(H)=\brk{12}$ and
$E_0=\{e_i:i\in \brk{7}\}$.
In the following, we will show that
$T^*$ is a tree whenever $T$ is a tree.

%\red{  the edges belonging to $E_0$ are drawn as    thick edges.}

\begin{figure}[h!]
\centering

\tikzp{0.95}
{
%%%%%%%%%%%%%%%%%%% for (c) F'.
\foreach \place/\z in
{{(0,0)/1},{(1,0)/4},{(2,0)/2},{(3,0)/8},
{(-1,1.5)/6},{(0,1.5)/7},{(1,1.5)/5},{(2,1.5)/10},
{(3,1.5)/3},{(4,1.5)/9}}
\node[cblue] (w\z) at \place {};

\filldraw[black] (w1) circle (0pt) node[anchor=north]
{$1$};
\filldraw[black] (w4) circle (0pt) node[anchor=north]
{\small $4$};
\filldraw[black] (w2) circle (0pt) node[anchor=north]
{\small $2$};
\filldraw[black] (w8) circle (0pt) node[anchor=north]
{\small $8$};

\filldraw[black] (w3) circle (0pt) node[anchor=south]
{\small $3$};
\filldraw[black] (w5) circle (0pt) node[anchor=south]
{\small $5$};
\filldraw[black] (w6) circle (0pt) node[anchor=south]
{\small $6$};
\filldraw[black] (w7) circle (0pt) node[anchor=south]
{\small $7$};
\filldraw[black] (w9) circle (0pt) node[anchor=south]
{\small $9$};
\filldraw[black] (w10) circle (0pt) node[anchor=south]
{\small $10$};

\draw[black] (w5) -- (w1);
\draw[black] (w4) -- (w6);
\draw[black] (w10) -- (w8) -- (w9);

\filldraw[] (1.5,-0.5) circle (0pt) node[anchor=north]
{\small (c) $F'\in \setf_{U,p'}$};

%\draw[thin]
%(w1) edge[ultra thick] node[above]{$e_1$} (w4)
%(w2) edge[ultra thick] node[above]{$e_2$} (w3);

%%%%%%%%%%%% for the upper part, (a) T

\foreach \place/\z in
{{(-1,3.5)/1},{(0,3.5)/7},{(1,3.5)/3},{(2,3.5)/4},
{(3,3.5)/8},{(4,3.5)/9},
{(0,5)/6},{(1,5)/5},{(2,5)/10},
{(3,5)/2},{(-2,4.25)/11},{(4.5,4.25)/12}}
\node[cblue] (v\z) at \place {};

\filldraw[black] (v1) circle (0pt) node[anchor=north]
{$1$};
\filldraw[black] (v4) circle (0pt) node[anchor=north]
{\small $4$};
\filldraw[black] (v2) circle (0pt) node[anchor=south]
{\small $2$};
\filldraw[black] (v8) circle (0pt) node[anchor=north]
{\small $8$};

\filldraw[black] (v3) circle (0pt) node[anchor=north]
{\small $3$};
\filldraw[black] (v5) circle (0pt) node[anchor=south]
{\small $5$};
\filldraw[black] (v6) circle (0pt) node[anchor=south]
{\small $6$};
\filldraw[black] (v7) circle (0pt) node[anchor=north]
{\small $7$};
\filldraw[black] (v9) circle (0pt) node[anchor=north]
{\small $9$};
\filldraw[black] (v10) circle (0pt) node[anchor=south]
{\small $10$};

\filldraw[black] (v11) circle (0pt) node[anchor=south]
{\small $11$};
\filldraw[black] (v12) circle (0pt) node[anchor=south]
{\small $12$};

\draw[black] (v5) -- (v1);
\draw[black] (v4) -- (v6);
\draw[black] (v8) -- (v10) -- (v9);

\draw[thin]
(v1) edge[ultra thick] node[xshift=2pt,left]{$e_1$} (v6)
(v2) edge[ultra thick] node[xshift=2pt,pos=0.6, left]{$e_2$} (v4)
(v3) edge[ultra thick] node[xshift=3pt, pos=0.4, left]{$e_3$} (v5)
(v6) edge[ultra thick] node[xshift=-3pt, pos=0.7,right]{$e_4$} (v7)
(v5) edge[ultra thick] node[yshift=-2pt, above]{$e_5$} (v10)
(v1) edge[ultra thick] node[xshift=2pt, left]{$e_6$} (v11)
(v9) edge[ultra thick] node[xshift=-1pt, right]{$e_7$} (v12);

\filldraw[] (1.5,3) circle (0pt) node[anchor=north]
{\small (a) $T=E_0\cup F$};

%%%%%%%%%%%%%%%%%%% (right hand side) figure (d) $T^*$

\foreach \place/\z in
{{(7,0)/1},{(8,0)/4},{(9,0)/2},{(10,0)/8},
{(6,1.5)/6},{(7,1.5)/7},{(8,1.5)/5},{(9,1.5)/10},
{(10,1.5)/3},{(11,1.5)/9},
{(5.5,0.75)/11},{(12,0.75)/12}}
\node[cblue] (w\z) at \place {};

\filldraw[black] (w1) circle (0pt) node[anchor=north]
{$1$};
\filldraw[black] (w4) circle (0pt) node[anchor=north]
{\small $4$};
\filldraw[black] (w2) circle (0pt) node[anchor=north]
{\small $2$};
\filldraw[black] (w8) circle (0pt) node[anchor=north]
{\small $8$};

\filldraw[black] (w3) circle (0pt) node[anchor=south]
{\small $3$};
\filldraw[black] (w5) circle (0pt) node[anchor=south]
{\small $5$};
\filldraw[black] (w6) circle (0pt) node[anchor=south]
{\small $6$};
\filldraw[black] (w7) circle (0pt) node[anchor=south]
{\small $7$};
\filldraw[black] (w9) circle (0pt) node[anchor=south]
{\small $9$};
\filldraw[black] (w10) circle (0pt) node[anchor=south]
{\small $10$};
\filldraw[black] (w11) circle (0pt) node[anchor=south]
{\small $11$};
\filldraw[black] (w12) circle (0pt) node[anchor=north]
{\small $12$};

\draw[black] (w5) -- (w1);
\draw[black] (w4) -- (w6);
\draw[black] (w10) -- (w8) -- (w9);

\filldraw[] (8.5,-0.5) circle (0pt) node[anchor=north]
{\small (d) $T^*=E_0\cup F'$};

\draw[thin]
(w1) edge[ultra thick] node[xshift=2pt,left]{$e_1$} (w6)
(w2) edge[ultra thick] node[yshift=-2pt, above]{$e_2$} (w4)
(w3) edge[bend left=25,ultra thick] node[yshift=2pt, below]{$e_3$} (w5)
(w6) edge[ultra thick] node[yshift=-2pt, above]{$e_4$} (w7)
(w5) edge[ultra thick] node[yshift=-2pt, above]{$e_5$} (w10)
(w1) edge[ultra thick] node[yshift=1pt, below]{$e_6$} (w11)
(w9) edge[ultra thick] node[xshift=4.5pt,yshift=-5.5pt, above]{$e_7$} (w12);

%\draw[thin]
%(w1) edge[ultra thick] node[above]{$e_1$} (w4)
%(w2) edge[ultra thick] node[above]{$e_2$} (w3);

%%%%%%%%%%%% for the upper part, (b) F

\foreach \place/\z in
{{(6,3.5)/1},{(7,3.5)/7},{(8,3.5)/3},{(9,3.5)/4},
{(10,3.5)/8},{(11,3.5)/9},
{(7,5)/6},{(8,5)/5},{(9,5)/10},{(10,5)/2}}
\node[cblue] (v\z) at \place {};

\filldraw[black] (v1) circle (0pt) node[anchor=north]
{$1$};
\filldraw[black] (v4) circle (0pt) node[anchor=north]
{\small $4$};
\filldraw[black] (v2) circle (0pt) node[anchor=south]
{\small $2$};
\filldraw[black] (v8) circle (0pt) node[anchor=north]
{\small $8$};

\filldraw[black] (v3) circle (0pt) node[anchor=north]
{\small $3$};
\filldraw[black] (v5) circle (0pt) node[anchor=south]
{\small $5$};
\filldraw[black] (v6) circle (0pt) node[anchor=south]
{\small $6$};
\filldraw[black] (v7) circle (0pt) node[anchor=north]
{\small $7$};
\filldraw[black] (v9) circle (0pt) node[anchor=north]
{\small $9$};
\filldraw[black] (v10) circle (0pt) node[anchor=south]
{\small $10$};

\draw[black] (v5) -- (v1);
\draw[black] (v4) -- (v6);
\draw[black] (v8) -- (v10) -- (v9);

\filldraw[] (8.5,3) circle (0pt) node[anchor=north]
{\small (b) $F\in \setf_{U,p}$};

}

\caption{$T$, $F$, $F'$ and $T^*$, where
 $E_0=\{e_i: i\in \brk{7}\}$
and $F\in \setf_{U,(6,4)}$}

\label{f17}
\end{figure}

\prop{pr4-1}
{ For any $E_0\subseteq E(H)$, if $H\spann{E_0}$ is an acyclic graph,
	then
\begin{enumerate}
\item
for any $T\in \sett_{U,p}(E_0)$,
	$T^*$ is a tree in $\sett_{U,p'}(E_0)$; and
\item
$T\rightarrow T^*$
is a bijection from $\sett_{U,p}(E_0)$ to $\sett_{U,p'}(E_0)$.
\end{enumerate}
}

\proof (i).  For $T=E_0\cup F$,
where $F\in \setf_{U,p}$,
by Lemma~\ref{le2-2} (iv),
$T=E_0\cup F$ is acyclic
if and only if $T^*=E_0\cup F'$ is acyclic.
Since
$V(T)=V(T^*)=V(H)$
and $|E(T)|=|E(T^*)|=|E_0|+|F|$,
(i) follows.
\\
(ii).
By Proposition~\ref{pr2-1},
$F\rightarrow F'$ is a bijection from
$\setf_{U,p}$ to $\setf_{U,p'}$,
implying that
$E_0\cup F\rightarrow E_0\cup F'$
is a bijection
from $\sett_{U,p}(E_0)$ to $\sett_{U,p'}(E_0)$
%by (\ref{set-p-E0}),
as desired, completing the proof.
\proofend

Recall that $\CT_{U,p}$  (defined in subsection 1.1) is the set of complete tiered graphs $Q$ in $\setg_{U,p}$ and
$\CT_{U,p}(H)=\{H\cup Q: Q\in \CT_{U,p}\}$.
Clearly, $\CT_{U,p}(H)$ has exactly ${|U|\choose p_1}$ graphs,  corresponding
to possible choices for the
subset of vertices $V_1(G)$ from $U$ for any graph $G\in \CT_{U,p}(H)$.     \label{r2-25}

For any two different graphs $G_1$ and $G_2$ in $\CT_{U,p}(H)$,
there may exist $T_i\in \ST(G_i)$ for $i=1,2$
such that $V(T_1)=V(T_2)$ and $E(T_1)=E(T_2)$.
For example, for $p=(2,2)$, if $H$ is the graph
with $V(H)=\brk{4}$ and $E(H)=\{e_1,e_2\}$,
where $e_i$ joins vertices
$i$ and $5-i$ for $i=1,2$,
then, $G_1$ and $G_2$ in Figure~\ref{f15} (a) and (c)
are graphs in $\CT_{U,p}(H)$, where $U=\brk{4}$ and $p=(2,2)$.
Observe that $G_i$ has a spanning $T_i$
with  $V(T_i)=\brk{4}$ and $E(T_i)=\{\{1,2\},e_1,e_2\}$
for $i=1,2$, as shown in Figure~\ref{f15} (b) and (d).
However, as $V_1(T_1)=\{2,3\}\ne \{2,4\}= V_1(T_2)$,
$T_1$ and $T_2$ are different quasi-tiered trees
by definition.

\begin{figure}[H]
\centering

\tikzp{1.1}
{
%%%%%%%%%%%%%%%%%%% figure (a)
\foreach \place/\z in
{{(0,0)/1},{(1.5,0)/4},
{(0,1.5)/2},{(1.5,1.5)/3}}
\node[cblue] (w\z) at \place {};

\filldraw[black] (w1) circle (0pt) node[anchor=north]
{$1$};
\filldraw[black] (w2) circle (0pt) node[anchor=south]
{\small $2$};
\filldraw[black] (w3) circle (0pt) node[anchor=south]
{\small $3$};
\filldraw[black] (w4) circle (0pt) node[anchor=north]
{\small $4$};

\draw[black] (w3) -- (w1) -- (w2);

\draw[thin]
(w1) edge[ultra thick] node[above]{$e_1$} (w4)
(w2) edge[ultra thick] node[above]{$e_2$} (w3);

%\filldraw[black] (2.5,-0.25) circle (0pt) node[anchor=north] {(a)};

%%%%%%%%%%%%%%%%%%% figure (b)
\foreach \place/\z in
{{(3.5,0)/1},{(5,0)/4},
{(3.5,1.5)/2},{(5,1.5)/3}}
\node[cblue] (w\z) at \place {};

\filldraw[black] (w1) circle (0pt) node[anchor=north]
{$1$};
\filldraw[black] (w2) circle (0pt) node[anchor=south]
{\small $2$};
\filldraw[black] (w3) circle (0pt) node[anchor=south]
{\small $3$};
\filldraw[black] (w4) circle (0pt) node[anchor=north]
{\small $4$};

\draw[black] (w1) -- (w2);

\draw[thin]
(w1) edge[ultra thick] node[above]{$e_1$} (w4)
(w2) edge[ultra thick] node[above]{$e_2$} (w3);

%%%%%%%%%%%%%%%%%%% figure (c)  +4.5
\foreach \place/\z in
{{(8,0)/1},{(9.5,0)/3},
{(8,1.5)/2},{(9.5,1.5)/4}}
\node[cblue] (w\z) at \place {};

\filldraw[black] (w1) circle (0pt) node[anchor=north]
{$1$};
\filldraw[black] (w3) circle (0pt) node[anchor=north]
{\small $3$};
\filldraw[black] (w2) circle (0pt) node[anchor=south]
{\small $2$};
\filldraw[black] (w4) circle (0pt) node[anchor=south]
{\small $4$};

\draw[black] (w2) -- (w1) -- (w4) -- (w3);

\draw[thin]
(w1) edge[ultra thick, bend right=20]
node[pos=0.3, below]{$e_1$} (w4)
(w2) edge[ultra thick] node[pos=0.3, above]{$e_2$} (w3);

%%%%%%%%%%%%%%%%%%% figure (d)  +3.5
\foreach \place/\z in
{{(11.5,0)/1},{(13,0)/3},
{(11.5,1.5)/2},{(13,1.5)/4}}
\node[cblue] (w\z) at \place {};

\filldraw[black] (w1) circle (0pt) node[anchor=north]
{$1$};
\filldraw[black] (w3) circle (0pt) node[anchor=north]
{\small $3$};
\filldraw[black] (w2) circle (0pt) node[anchor=south]
{\small $2$};
\filldraw[black] (w4) circle (0pt) node[anchor=south]
{\small $4$};

\draw[black] (w2) -- (w1); % -- (w4) -- (w3);

\draw[thin]
(w1) edge[ultra thick, bend right=20]
node[pos=0.3, below]{$e_1$} (w4)
(w2) edge[ultra thick] node[pos=0.3, above]{$e_2$} (w3);

}
 {}\hspace{0.5 cm} (a) $G_1$ \hspace{1.6 cm}
 (b) $T_1\in \ST(G_1)$ \hspace{2.7 cm}
 (c) $G_2$ \hspace{1.6 cm}
 (d) $T_2\in \ST(G_2)$

\caption{$G_1,G_2\in \CT_{U, (2,2)}(H)$
with $T_i\in \ST(G_i)$
for $i=1,2$ such that $E(T_1)=E(T_2)$,
where  $V(H)=U=\brk{4}$  and the thick edges belong to $H$}
\label{f15}
\end{figure}

\lemma{le4-4}
{
Let $H$ be a graph
%(possibly not tiered)
 with $U\subseteq V(H)$.
For any two different graphs $G_1$ and $G_2$ in $\CT_{U,p}(H)$, we have
$\ST(G_1)\cap \ST(G_2)=\emptyset$.
}

\proof Suppose that $G_i=H\cup Q_i$, where $Q_i\in \CT_{U,p}$ for $i=1,2$.
As $G_1$ and $G_2$ are different graphs in $\CT_{U,p}(H)$,
we have $V_1(Q_1)\ne V_1(Q_2)$.
Thus, if $T_i\in \ST(G_i)$ for $i=1,2$, then
$V_1(T_1)=V_1(Q_1)\ne V_1(Q_2)=V_1(T_2)$ holds.
This completes the proof.
 \proofend

By Proposition~\ref{pr4-1} (ii), the next consequence follows.

\prop{pr4-2}
{
For any graph $H$
 with $U\subseteq V(H)$,
$T\rightarrow T^*$ is a bijection
from $\bigcup\limits_{G\in \CT_{U,p}(H)}\ST(G)$
to $\bigcup\limits_{G\in \CT_{U,p'}(H)}\ST(G)$.
}

\proof For any $E_0\subseteq E(H)$,
by Proposition~\ref{pr4-1} (ii),
$T\rightarrow T^*$
for $T\in \sett_{U,p}(E_0)$
is a bijection from $\sett_{U,p}(E_0)$ to $\sett_{U,p'}(E_0)$.
Note that
$$
\bigcup_{G\in \CT_{U,p}(H)}\ST(G)
=\bigcup_{E_0\subseteq E(H)}\sett_{U,p}(E_0),
$$
where the equality holds if $p$ is replaced by $p'$.
Thus, the result follows.
\proofend

 \begin{rem}
 	For an ordered partition
$p=(p_1,\cdots,p_{r},p_{r+1}, \cdots,p_m)$ of $n$,  let
 $\pi_r(p)$ denote the ordered partition $(p_1,\cdots,p_{r-1},p_{r+1},p_r, p_{r+2},\cdots,p_m)$ of $n$ obtained by exchanging $p_r$ and $p_{r+1}$.
  When restricted to $\sett_p$,
  the bijection $T\rightarrow T^*$  recovers the  bijection $\phi_r:
  \sett_p  \rightarrow \sett_{\pi_r(p)}$
  which was established by Yan et al. \cite{yan2020bijective}.
 \end{rem}

\section{Preliminaries on external activities in quasi-tiered graphs
%Proof of a special case of Theorem~\ref{main-th3}
%An identity on Tutte polynomials
\label{sec5}
}

\subsection{A property on external activity of a spanning tree
	\label{sec5-01}
}

In this subsection, all graphs considered are not necessarily tiered.
The Tutte polynomial $\tutte_G(x,y)$ of a graph $G$
has many different expressions
\cite{bjorner1992homology,brylawski1992tutte, crapo1969tutte,
	tutte1954contribution, welsh1999tutte}.
Recall that
$\ST(G)$ denotes the set of spanning trees of $G$.  Let $\omega$ be an injective function from $E(G)$ to the set $\R$ of real numbers.   \label{r2-6}
 Clearly, $\omega$ defines a total order on the edges of $G$.
If $G$ is connected,
$\tutte_G(x,y)$ can be expressed in terms
of spanning trees~\cite{crapo1969tutte, tutte1954contribution}:
\begin{equation}\label{eq1-2}
	\tutte_G(x,y)=\sum_{T\in \ST(G)}x^{ia(T)}y^{ea(T)},
\end{equation}
where $ia(T)$ and $ea(T)$, introduced below,
are respectively the internal and external activities
of $T$ with respect to an
injective weight function $\omega$ on $E(G)$.

Assume that $G$ is connected.
For any subgraph $H$ of $G$ and
any edge $e\in E(G)\setminus E(H)$,
let $H\cup e$ denote the subgraph
of $G$ with vertex set $V(H)\cup \{u,v\}$
and edge set $E(H)\cup \{e\}$,
where $u$ and $v$ are the endpoints of $e$.
If $H$ is a forest, then $H\cup e$
contains a cycle \iff the endpoints of $e$ are in the same component of $F$.

	The definitions of $ia(T)$ and $ea(T)$ for a spanning tree $T$ of $G$ can be
	extended to $ia(F)$ and $ea(F)$ for a forest $F$ of $G$ respectively.
	%For a spanning forest $F$ of $G$,
	An edge $e$ in $F$ is called
	{\it internally active} with respect to $(\omega, F)$ (or simply $F$)
	if $\omega(e)<\omega(e')$ holds
	for each $e'\in E(G)\setminus E(F)$
	with the properties that
	$e$ is on a cycle of $F\cup e'$,
	and
	an edge $e\in E(G)\setminus E(F)$ is said to be
	{\it externally active} with respect to
	$(\omega, F)$ (or simply $F$)
	if $F\cup e$ contains a cycle $C$ and
	$\omega(e)<\omega(e')$ holds for all edges $e'$ in $E(C)\setminus \{e\}$.
	Let $ia(F)$ (resp. $ea(F)$) denote the number of
	internally (resp. externally) active edges
	with respect to $F$.
	We also call $ia(F)$ and $ea(F)$ the {\it internal activity} and
	{\it external activity} of $F$ in $G$ respectively.

\iffalse
For any $E_0\in E(G)$, if $G\spann{E_0}$
is a forest, we also write
$ia(E_0)$ and $ea(E_0)$ for
$ia(F)$ and $ea(F)$ respectively.
\fi
Note that, for $T\in \ST(G)$,
$ia(T)$ and $ea(T)$ are associated
with the weight function $\omega$ and $G$,
although $\tutte_G(x,y)$ is independent of the weight function $\omega$.
Thus, $ia(T)$ and $ea(T)$ are respectively written as
$ia_{\omega,G}(T)$ and $ea_{\omega,G}(T)$ whenever there is a danger of
confusion. \label{r1-3}

	For any two graphs $G$ and $F$ with  $V(F)\subseteq V(G)$, denote by $G\bullet F$ the graph obtained from $G$
	by removing all edges in $ E(F)\cap E(G)$ %from $G$
	and identifying all vertices in $F_i$
	for all $i=1,2,\cdots,k$,
	where $F_1, F_2,\cdots,F_k$ are the components of $F$.
	Thus, $G\bullet F$ is a graph with $k$ vertices when $V(F)= V(G)$.
For example, if $G$ is the cycle graph $C_n$
	and $F$ is a spanning subgraph of $G$ with $k$ components, then
	then $G\bullet F\cong C_{k}$.
Note that in the definition of $G\bullet F$,
it is possible that $E(F)\not\subseteq E(G)$.

For any partition $(E_1,E_2)$ of $E(G)$,
when  $\omega(e_1)<\omega(e_2)$ holds
for all $e_1\in E_1$ and $e_2\in E_2$,
we have the following expression
for $ea_{\omega,G}(T)$
%in terms of$ea_{\omega,G_1\bullet F }(T\bullet F)$ and $ea_{\omega,G_2}(F)$, where

%A property on $ea_{\omega,G}(T)$ is presented below.

\prop{pr02-1}
{Let $G$ be a connected graph with an injective weight function
	$\omega$ on $E(G)$ and $T\in \ST(G)$.
	 If $E(G)$ is partitioned into $E_1$ and $E_2$
		such that
		$\omega(e_1)<\omega(e_2)$ holds for all $e_1\in E_1$ and
		$e_2\in E_2$,
	then
	\begin{equation}\label{eq02-2}
		ea_{\omega,G}(T)=
		ea_{\omega,G_1\bullet F }(T\bullet F)
		+ea_{\omega,G_2}(F),
	\end{equation}
	where $G_i=G\langle E_i\rangle$ and $F=T\spann{E(T)\cap E_2}$.
}

\proof For any $e\in E(G)\setminus E(T)$, let
$C_e$ be the unique cycle in $T\cup e$.
Let $\EX_G(T)$ denote the set of edges in $E(G)\setminus E(T)$
which are externally active with respect  to $T$.
We prove this result by showing the following claims.

\inclaim For $e\in E_1\setminus E(T)$,
$e\in \EX_G(T)$ \iff $e\in \EX_{G_1\bullet F}(T\bullet  F)$.

As $E(F)=E(T)\cap E_2$,
$T\bullet F$ is a spanning tree of
$G_1\bullet  F$ and $C'_e$ is the unique cycle
of $(T\bullet F)\cup e$,
where $C_e'$ is the graph $C_e\slash (E_2\cap E(C_e))$.
 Then  each pair of consecutive statements below are equivalent:\\
	(a). $e\in \EX_G(T)$;  \\
	(b). $\omega(e)\le \omega(e')$ holds for all $e'\in E(C_e)$; \\
	(c). $\omega(e)\le \omega(e')$ holds for all $e'\in E(C_e')$; and \\
	(d).  $e\in \EX_{G_1 \bullet F}(T\bullet F)$. \\
	Thus, Claim 1 holds.

\inclaim For any $e\in E_2\setminus E(T)$,
$e\in \EX_G(T)$ \iff
$e\in \EX_{G_2}(F)$.

 First, by the given condition $\omega(e)>\omega(e')$ for all $e'\in E_1$,
	one can  easily verify  that each pair of
	consecutive statements below are equivalent:\\
	(a'). $e\in \EX_G(T)$;\\
	(b'). $\omega(e)\le \omega(e')$ holds for all
	$e'\in E(C_e)$ and $E(C_e)\setminus \{e\}
	\subseteq E_2\cap E(T)=E(F)$; and \\
	(c'). $e\in \EX_{G_2}(F)$.\\
	This completes the proof of Claim 2. Combining  Claims 1 and 2, we are led to
	(\ref{eq02-2}).
\proofend

For any graph $H$ and any proper subset $S$ of $V(H)$, let $H-S$ be the graph obtained from $H$ by deleting all vertices in $S$ and all edges incident with vertices in $S$.
Note that if $F$ is a subgraph of $H$
and $S\subsetneq V(F)$ is a set of isolated vertices of $F$, then
$H\bullet F$ and $H\bullet (F-S)$
are the same graph.
Thus, the following conclusion holds.

\begin{cor}\label{co3-1}
Proposition~\ref{pr02-1} still holds
if $G_2$ and $F$ are  replaced by $G\spann{E_2}-S$ and  $T\spann{E(T)\cap E_2}-S$ respectively, where
$S\subsetneq V(G)$ is a set of isolated vertices in
$G\spann{E_2}$.
\end{cor}

%Proposition~\ref{pr02-1} can be extended to a result on matroids \cite{oxley2006matroid}.

\subsection{Identities on external activities in quasi-tiered graphs
	\label{sec5-2}
}

 In the  remainder of this section,
we assume that $U$ is a non-empty finite set of positive integers,
$p=(p_1,p_2)$ is an ordered partition $p=(p_1,p_2)$ of $|U|$,
 $p'=(p_2,p_1)$  and
  $H$ is a graph with $U\subseteq V(H)$.
%and $E(H)=\{e_i:i\in \brk{m}\}$.

We first define two weight functions
$\omega_1$ and $\omega_2$.
Let $\omega_1$ be an injective weight function on the set
$E(H)\cup \bigcup_{G\in \CT_{U,p}}E(G)$ defined below:  \label{r2-27}
\begin{enumerate}
\item  $\omega_1(e)<0$ for each $e\in E(H)$,
and $\omega_1(e_1)\ne \omega_1(e_2)$ for any different edges $e_1,e_2$ in $E(H)$; and
\item for any edge $e$ joining vertex $u$ and vertex $v$, where $1\le u<v$,
if $e\notin E(H)$, then $\omega_1(e)=u+\frac v{N}$, where $N=\max\{w:w\in U\}$.\footnote{$\omega_1$ is actually the lexicographic order on $(\min \{u,v\}, \max\{u,v\})$) for edges $e=uv$ which don't belong to $E(H)$.}
\end{enumerate}
Let $\omega_2$ be the weight function on the set
$E(H)\cup \bigcup_{G\in \CT_{U,p'}}E(G)$
defined as follows:
\begin{enumerate}
\item $\omega_2(e)=\omega_1(e)$ for each $e\in E(H)$; and
\item for any edge $e$ joining vertex $u$ and vertex $v$, where $1\le u<v$,
if $e\notin E(H)$, then $\omega_2(e)=(N+1-v)+\frac {N+1-u}N$.
\footnote{$\omega_2$ is actually the lexicographic order on $(-\max \{u,v\}, -\min\{u,v\})$) for edges $e=uv$ which don't belong to $E(H)$.}
\end{enumerate}

By definitions, it is easily seen that
$\omega_1$ and $\omega_2$
have the following properties.

\lemma{le5-1}
{Let $Q\in \CT_{U,p}$ and $E^*\subseteq E(Q)$.
For $e=uv\in E^*$, where $u<v$,
$\omega_1(e)=\min\limits_{e'\in E^*}\omega_1(e')$
\iff $u=\min \{u'\in V(Q): u'
\mbox{ is an end of some edges in } E^*\}$
and $v=\min \{v'\in V(Q): uv' \in E^*\}$.
}

\lemma{le5-2}
{Let $Q\in \CT_{U,p'}$ and $E^*\subseteq E(Q)$.
For $e=uv\in E^*$, where $u<v$,
$\omega_2(e)=\min\limits_{e'\in E^*}\omega_2(e')$
\iff $v=\max \{v'\in V(Q): v' \mbox{ is an end of some edges in } E^*\}$
and $u=\max \{u'\in V(Q):  vu' \in E^*\}$.
}

Recall that, in Section~\ref{sec3},
we introduce the quasi-tiered tree $T^*$
for any tree $T=E_0\cup F$, where $E_0\subseteq E(H)$ and $F\in \setf_{U,p}$.
By Proposition~\ref{pr4-1} (ii),
$T\rightarrow T^*$ is a bijection from $\sett_{U,p}(E_0)$ to $\sett_{U,p'}(E_0)$.
In this section, our main purpose is to show that $T$ and $T^*$ have the same external activity with respect the $\omega_1$ and $\omega_2$ respectively.

We first establish a similar result for
a forest $F\in \setf_{U,p}$ and its dual
forest $F'\in \setf_{U,p'}$.
For any graph $Q\in \setg_{U,p}$,
there exists a unique complete tiered graph in $\CT_{U,p}$, denoted by $CT(Q)$,
with $V_i(CT(Q))=V_i(Q)$ for $i=1,2$.
We say that $CT(Q)$ is the {\it complete tiered graph determined} by $Q$.

\prop{pr5-1}
{For any forest $F\in \setf_{U,p}$,
$
ea_{\omega_1,CT(F)}(F)=ea_{\omega_2,CT(F')}(F').
$
}

\proof Let $F_1,F_2,\cdots,F_k$ be the components of $F$.
By the definition of $F'$   \label{r2-7}
and Lemma~\ref{le2-1} (i),
 $F'$  has components $F'_1,F'_2,\cdots, F'_k$.
By the definition of external activity,
\begin{equation}\label{pr5-1-e2}
ea_{\omega_1,CT(F)}(F)=\sum_{i=1}^k  ea_{\omega_1,CT(F)}(F_i),
\qquad
ea_{\omega_2,CT(F')}(F')=\sum_{i=1}^k  ea_{\omega_2,CT(F')}(F'_i).
\end{equation}
By (\ref{pr5-1-e2}), to prove %(\ref{pr5-1-e1}),
 the assertion of Proposition~\ref{pr5-1},
it suffices to show that for each $i=1,2,\cdots,k$,
$ea_{\omega_1,CT(F)}(F_i)=ea_{\omega_2,CT(F')}(F'_i)$ holds.

Without loss of generality, we will show that
$ea_{\omega_1,CT(F)}(F_1)=ea_{\omega_2,CT(F')}(F'_1)$.
By definition,
\iffalse
$ea_{\omega_1,CT(F)}(F_1)$ is the number of edges
$e$ in $E(CT(F))\setminus E(F_1)$
such that $F_1\cup e$ has a cycle $C$
and $\omega_1(e)\le \omega_1(e')$ holds for all $e'\in E(C)$.
Thus,
\fi
to compute $ea_{\omega_1,CT(F)}(F_1)$,
we need only to consider those edges in
$E(CT(F))\setminus E(F_1)$
with both ends in $V(F_1)$.

Assume that $V(F_1)=\{x_1,x_2,\cdots,x_r\}$,
where $1\le x_1<x_2<\cdots<x_r$.

  \inclaim For $1\le i<j\le r$,
$x_{i}x_{j}\in  E(CT(F))\setminus E(F_1)$
\iff $x_{r+1-i}x_{r+1-j}\in  E(CT(F'))\setminus E(F'_1)$.\\
Claim 1 follows directly from the construction of $F'$ and $CT(F')$.

\inclaim For  $1\le i<j\le r$,
if $x_{i}x_{j}\in E(CT(F))\setminus E(F_1)$,
then $x_{i}x_{j}$ is externally active with respect to
$(\omega_1,F_1)$
\iff $x_{r+1-i}x_{r+1-j}$ is externally active with respect to
$(\omega_2, F'_1)$.

Assume that $x_{i}x_{j}\in E(CT(F))\setminus E(F_{1})$ with $x_i,x_j\in V(F_1)$. By Claim 1,
$x_{r+1-i}x_{r+1-j}\in E(CT(F'))\setminus E(F'_{1})$, and by the definition of $F_1'$,
$x_{r+1-i},x_{r+1-j}\in V(F')$.

Let $e=x_ix_j$ and let $P: x_{z_1}x_{z_2}\cdots x_{z_w}$ be
the unique path in $F_1$ joining $x_{i}$ and $x_{j}$,
where $z_1=i$, $z_w=j$
and $z_1,z_2,\cdots,z_w$ are numbers in $\brk{r}$.
Then, the unique cycle $C$ in $F_1\cup e$ is
the cycle $x_{z_1}x_{z_2}\cdots x_{z_w} x_{z_1}$.
By Lemma~\ref{le5-1},
$e$ is externally active with respect to
$(\omega_1, F_1)$ \iff
$x_i<\min\{x_{z_s}:2\le s\le w\}$
and $x_{z_w}<x_{z_2}$,
i.e., $i<\min\{z_s:2\le s\le w\}$ and $z_w<z_2$.
Thus, the following subclaim holds.

\noindent {\bf Claim 2.1}:
$x_{i}x_{j}$ is externally active respect to
$(\omega_1,F_1)$ \iff $i<\min\{z_s:2\le s\le w\}$
and $z_w<z_2$.

By the definition of $F'_1$ and
the fact that $P$ is a path in $F_1$,
$x_{r+1-z_1}x_{r+1-z_2}\cdots x_{r+1-z_w}$ is a path in $F'_1$,
denoted by $P'$.
As $x_ix_j\in E(CT(F))\setminus E(F_1)$,
by  Claim 1, $x_{r+1-i}x_{r+1-j}\in E(CT(F'))\setminus E(F'_1)$.
Thus, $x_{r+1-z_1}x_{r+1-z_2}\cdots x_{r+1-z_w}x_{r+1-z_1}$
is the unique cycle in
$F'_1\cup x_{r+1-z_1}x_{r+1-z_w}$,  denoted by $C'$,
as shown in Figure~\ref{f13} (b).
By Lemma~\ref{le5-2},
$x_{r+1-i}x_{r+1-j}$ is externally active respect to
$(\omega_2,F'_1)$ \iff
$r+1-i>\max\{r+1-z_s: s=2,3,\cdots,w\}$
and $r+1-z_w>r+1-z_2$,
implying the following subclaim.

\noindent {\bf Claim 2.2}:
$x_{r+1-i}x_{r+1-j}$ is externally active respect to
$(\omega_2,F'_1)$ \iff $i<\min\{z_s:2\le s\le w\}$
and $z_w<z_2$.

Claim 2 follows from Subclaims 2.1 and 2.2.

By Claims 1 and 2,
$ea_{\omega_1,CT(F)}(F_1)=ea_{\omega_2,CT(F')}(F'_1)$ holds.
Hence the result holds.
\proofend

\begin{figure}[h!]
\centering

\tikzp{1}
{
%%%%%%%%%%%%%%%%%%% figure (a)
\foreach \place/\z in
{{(0,0)/1},{(1.5,0)/2},{(3,0)/3},{(5,0)/4},
{(0.75,2)/5},{(2.25,2)/6},{(3.75,2)/7},{(5.75,2)/8}}
\node[cblue] (w\z) at \place {};

\filldraw[black] (w1) circle (0pt) node[anchor=north]
{$x_{z_1}$};
\filldraw[black] (w2) circle (0pt) node[anchor=north]
{$x_{z_3}$};
\filldraw[black] (w3) circle (0pt) node[anchor=north]
{$x_{z_5}$};
\filldraw[black] (w4) circle (0pt) node[anchor=north]
{$x_{z_{w-1}}$};
\filldraw[black] (w5) circle (0pt) node[anchor=south]
{$x_{z_2}$};
\filldraw[black] (w6) circle (0pt) node[anchor=south]
{$x_{z_4}$};
\filldraw[black] (w7) circle (0pt) node[anchor=south]
{$x_{z_6}$};
\filldraw[black] (w8) circle (0pt) node[anchor=south]
{$x_{z_w}$};

\draw[black] (w1) -- (w5) -- (w2) -- (w6) -- (w3) -- (w7);
\draw[black] (w8) -- (w4);
\draw[dashed] (w4) -- (w7);

\draw[thin]
(w1) edge[thick] node[above]{$e$} (w8);

\filldraw[black] (4,0) circle (0pt) node[]
{$\cdots$};
\filldraw[black] (4.75, 2) circle (0pt) node[]
{$\cdots$};

\iffalse
\draw[thin]
(w1) edge[ultra thick] node[above]{$e_1$} (w4)
(w2) edge[ultra thick] node[above]{$e_2$} (w3);
\fi

%\filldraw[black] (2.5,-0.25) circle (0pt) node[anchor=north] {(a)};

%%%%%%%%%%%%%%%%%%% figure (b)
\foreach \place/\z in
{{(9.5,0)/1},{(11,0)/2},{(12.5,0)/3},{(14.5,0)/4},
{(8.75,2)/5},{(10.25,2)/6},{(11.75,2)/7},{(13.75,2)/8}}
\node[cblue] (w\z) at \place {};

\filldraw[black] (w1) circle (0pt) node[anchor=north]
{$x_{r+1-z_2}$};
\filldraw[black] (w2) circle (0pt) node[anchor=north]
{$x_{r+1-z_4}$};
\filldraw[black] (w3) circle (0pt) node[anchor=north]
{$x_{r+1-z_6}$};
\filldraw[black] (w4) circle (0pt) node[anchor=north]
{$x_{r+1-z_w}$};
\filldraw[black] (w5) circle (0pt) node[anchor=south]
{$x_{r+1-z_1}$};
\filldraw[black] (w6) circle (0pt) node[anchor=south]
{$x_{r+1-z_3}$};
\filldraw[black] (w7) circle (0pt) node[anchor=south]
{$x_{r+1-z_5}$};
\filldraw[black] (w8) circle (0pt) node[anchor=south]
{$x_{r+1-z_{w-1}}$};

\draw[black] (w3) -- (w7) -- (w2) -- (w6) -- (w1) -- (w5);
\draw[dashed] (w3) -- (w8);
\draw[black] (w4) -- (w8);

\draw[thin]
(w4) edge[thick] node[below]{$e'$} (w5);'

\filldraw[black] (13.5,0) circle (0pt) node[]
{$\cdots$};
\filldraw[black] (12.75, 2) circle (0pt) node[]
{$\cdots$};

}
 (a) Cycle $C$ in $F_1\cup e$
\hspace{4.6 cm}
(b) Cycle $C'$ in $F_1'\cup e'$

\caption{Cycles $C$ in $F_1\cup e$  and $C'$ in $F'_1\cup e'$,
where $e=x_{z_1}x_{z_w}$ and $e'=x_{r+1-z_1}x_{r+1-z_w}$}
\label{f13}
\end{figure}

 For any tree $T\in \sett_{U,p}(E_0)$, where $E_0\subseteq E(H)$,
we have $T=E_0\cup F$ for some
$F\in \setf_{U,p}$.
\iffalse
By the definition of dual quasi-tired trees in Section~\ref{sec3},
$T^*$ is the graph $E_0\cup F'$,
where $F'\in \setf_{U,p'}$ is the dual forest of $F$.
\fi
By Proposition~\ref{pr4-1} (i),
$T^*=E_0\cup F'\in \sett_{U,p'}(E_0)$.
By Lemma~\ref{le4-4}, there is a unique graph in $\CT_{U,p'}(H)$,
denoted by $G^*_T$,
such that $T^*\in \ST(G^*_T)$.
Note that $H\cup CT(F')$ is a graph in the set $\CT_{U,p'}(H)$ and $T^*$ is a spanning tree
of $H\cup CT(F')$.
Thus,  $G^*_T$ is the graph $H\cup CT(F')$.

Observe that $G^*_T$ can be determined by the following procedure:
\begin{align*}
	T
	&\quad
	\Rightarrow
	\quad
	F=(U, E(T)\setminus E_{0}) \in {\mathcal{F}}_{U,p}\\
	&\quad
	\Rightarrow
	\quad
	F'\in {\mathcal{F}}_{U,p'}
	\\
	&\quad
	\Rightarrow
	\quad
	CT(F')\in {\mathcal{CG}}_{U,p'}\\
	&\quad
	\Rightarrow
	\quad
	G^{*}_{T}=H\cup CT(F')
\end{align*}
By applying Propositions~\ref{pr02-1} and~\ref{pr5-1}, we can now show that
for  $T\in \ST(G)$, where $G\in \CT_{U,p}(H)$,
$T$ and $T^*$ always have the same external activity with respect to $\omega_1$
and $\omega_2$ respectively.

\prop{pr5-2}
{For any $G\in \CT_{U,p}(H)$ and $T\in \ST(G)$,
$
ea_{\omega_1,G}(T)=ea_{\omega_2,G^*_T}(T^*).
$
}

\proof Assume that $G=H\cup Q$, where
$Q\in \CT_{U,p}$.
For $i=1,2$,
by the definition of $\omega_i$,
$\omega_i(e_1)<\omega_i(e_2)$ holds
for all $e_1\in E(H)$ and $e_2\in E(Q)$.
 Let $T\in \ST(G)$.
	Then $T=E_0\cup F$ for some $E_0\subseteq E(H)$
	 and $F=(U,E(T)\setminus E_0)\in \setf_{U,p}$.
Observe that  $Q=G\spann{E(Q)}-S$ and $F=T\spann{E(Q)\cap E(T)}-S$, where $S=V(H)\setminus U$ is a set of isolated vertices in $G\spann{E(Q)}$.
Then, by Corollary~\ref{co3-1},  we deduce that
\begin{equation}\label{pr5-2-e1}
ea_{\omega_1,G}(T)=
ea_{\omega_1,H\bullet F}(T\bullet F)
+ea_{\omega_1,Q}(F).
%+\sum_{1\le i\le r}|E(H[V(F_i)])|-\ell(H),
\end{equation}
 Note that $G^*_T=H\cup Q'$ and  $T^*=E_0\cup F'$, where $Q'=CT(F')$ and $F'$ is the dual forest of $F$.  Moreover, we have $Q'=G^*_T\spann{E(Q')}-S$ and $F'=T^*\spann{E(Q')\cap E(T^*)}-S$, where $S=V(H)\setminus U$ is a set of isolated vertices in $G^*_T\spann{E(Q')}$.
By  Corollary~\ref{co3-1} again,
%Proposition~\ref{pr02-1}
\begin{equation}\label{pr5-2-e2}
ea_{\omega_2,G^*_T}(T^*)=
ea_{\omega_2,H\bullet F'}(T^*\bullet F')
+ea_{\omega_2,Q'}(F').
%+\sum_{1\le i\le r}|E(H[V(F_i)])|-\ell(H),
\end{equation}

By Proposition~\ref{pr5-1}, we have
%\begin{equation}\label{eq2-n}
$ea_{\omega_1,Q}(F)=ea_{\omega_2,Q'}(F')$.
%\end{equation}
By (\ref{pr5-2-e1}) and (\ref{pr5-2-e2}),
it suffices to show that
\begin{equation}\label{eq1-n}
	ea_{\omega_1,H\bullet  F }(T\bullet F)
	=ea_{\omega_2,H\bullet  F' }(T^*\bullet F').
\end{equation}

Assume that $F_1,F_2,\cdots, F_k$ are the components of $F$.
Then, by definition and Lemma~\ref{le2-1} (i),
$F$ has components $F'_1,F'_2,\cdots, F'_k$.

By the definition of $F'$  and Lemma \ref{le2-2} (iv),
$V(F_i)=V(F'_i)$ for all $i=1,2,\cdots,k$,
implying that
$\sets(F)$ and $\sets(F')$ are the same partition of $U$, where
$\sets(F)=\{V(F_i): i=1,2,\cdots,k\}$
and $\sets(F')=\{V(F'_i): i=1,2,\cdots,k\}$.
Thus, $H\bullet F$ and $H\bullet F'$ are the same graph.

%Let $E_0=E(H)\cap E(T)$.
Obviously, $E_0=E(H)\cap E(T^*)$.
Observe that $T\bullet F=H\langle E_0\rangle \bullet F$ and  $T^*\bullet F'=H\langle E_0\rangle \bullet F'$.
Since $\sets(F)=\sets(F')$,
$H\langle E_0\rangle \bullet  F $
and $H\langle E_0\rangle \bullet  F' $
are the same graph,
implying that $T\bullet F$ and $T^*\bullet F'$ are the same tree.
By the definition of $\omega_1$ and $\omega_2$,
$\omega_1(e)=\omega_2(e)$ holds for all $e\in E(H)$. It follows that  (\ref{eq1-n}) holds.
Therefore,
$ea_{\omega_1,G}(T)=ea_{\omega_2,G^*_T}(T^*)
$ as desired.
\proofend

\section
{An elementary proof of Theorem~\ref{main-th3} via Tutte polynomials
\label{sec6}
}

\label{r2-28}

Let $U$ be a non-empty finite set of positive integers, let $p=(p_1,p_2)$ be an ordered partition of $|U|$,
$p'=(p_2,p_1)$
and let $H$ be a graph
with $U\subseteq V(H)$.

For $E_0\subseteq E(H)$,
let $\Phi_{U,p}(H,E_0)$ denote the set of ordered pairs
$(G,T)$, where $G\in \CT_{U,p}(H)$ and $T\in \ST(G)\cap \sett_{U,p}(E_0)$.
Clearly, $\Phi_{U,p}(H,E_0)=\emptyset$ if $H\langle E_0\rangle$ is not acyclic.
By Lemma~\ref{le4-4}, for any $T\in \sett_{U,p}(E_0)$,
there is only one graph $G\in \CT_{U,p}(H)$ such that
$T\in \ST(G)$ and $E(T)\cap E(H)=E_0$.
Thus, $\Phi_{U,p}(H,E_0)$ and $\sett_{U,p}(E_0)$ have the
same cardinality.

\prop{pr5-3}
{
%Let $H$ be a graph with $V(H)=U$.
For any $E_0\subseteq E(H)$,
\begin{equation}\label{eq5-5}
\sum_{(G,T)\in \Phi_{U,p}(H,E_0)}y^{ea_{\omega_1,G}(T)}
=\sum_{(G,T)\in \Phi_{U,p'}(H,E_0)} y^{ea_{\omega_2,G}(T)}.
\end{equation}
}

\proof The result is trivial when $H\langle E_0\rangle$
contains cycles, as both $\Phi_{U,p}(H,E_0)$ and $\Phi_{U,p'}(H,E_0)$
are empty in this case. Now assume that $H\langle E_0\rangle$ is acyclic.

By Proposition~\ref{pr4-1} (ii), $T\rightarrow T^*$
is a bijection from $\sett_{U,p}(E_0)$
to $\sett_{U,p'}(E_0)$,
implying that $(G,T)\rightarrow (G^*_T,T^*)$ is a
bijection from $\Phi_{U,p}(H,E_0)$ to $\Phi_{U,p'}(H,E_0)$.
By Proposition~\ref{pr5-2},
$ea_{\omega_1,G}(T)=ea_{\omega_2,G^*_T}(T^*)$ holds
for each $(G,T)\in \Phi_{U,p}(H,E_0)$. Hence,
\begin{equation}\label{eq5-4}
\sum_{(G,T)\in \Phi_{U,p}(H,E_0)} y^{ea_{\omega_1,G}(T)}
=\sum_{(G^{*}_T,T^{*})\in \Phi_{U,p'}(H,E_0)} y^{ea_{\omega_2,G^*_T}(T^{*})}
=\sum_{(G,T)\in \Phi_{U,p'}(H,E_0)} y^{ea_{\omega_2,G}(T)},
\end{equation}
by which (\ref{eq5-5}) follows.
\proofend

We are now going to prove Theorem~\ref{main-th3}.

\vspace{0.3 cm}

\noindent {\it Proof of Theorem~\ref{main-th3}}:
By (\ref{eq1-2}),
\begin{equation}\label{eq5-6}
	\sum_{G\in \CT^c_{U,p}(H)} \tutte_G(1,y)=
	\sum_{G\in \CT^c_{U,p}(H)}
	\sum_{T\in \ST(G)} y^{ea_{\omega_1,G}(T)}
	=\sum_{E_0\subseteq E(H)}
	\sum_{(G,T)\in \Phi_{U,p}(H,E_0)} y^{ea_{\omega_1,G}(T)};
\end{equation}
\begin{equation}\label{eq5-7}
	\sum_{G\in \CT^c_{U,p'}(H)}
	\tutte_G(1,y)=
	\sum_{G\in \CT^c_{U,p'}(H)} \sum_{T\in \ST(G)} y^{ea_{\omega_2,G}(T)}
	=\sum_{E_0\subseteq E(H)}
	\sum_{(G,T)\in \Phi_{U,p'}(H,E_0)} y^{ea_{\omega_2,G}(T)}.
\end{equation}
Observe that (\ref{eq1-6}) follows directly from
(\ref{eq5-6}), (\ref{eq5-7}) and Proposition~\ref{pr5-3}.
\proofend

By Corollary~\ref{main-cor}, Theorem~\ref{main-th} follows from
Theorem~\ref{main-th3}.

\section{An elementary proof of  Theorem~\ref{main-th4} via Tutte polynomials
\label{sec7}
}

\label{r2-31}

In this section, an elementary proof of Theorem~\ref{main-th4}
by applying Tutte polynomial
is provided.

For pairwise disjoint non-empty finite sets $U_1,U_2,\cdots,U_k$ of positive integers,
where $k\ge 2$,
let $CT(U_1,\cdots,U_k)$ denote
the complete tiered graph $G$ with $k$ tiers and
$V_i(G)=U_i$ for all $i=1,2,\cdots,k$.

For a graph $G$,
let $\tutte^c_G(y)=\tutte_G(1,y)$ when $G$ is connected,
and $\tutte^c_G(y)=0$ otherwise.
The following two lemmas will play essential roles in the proof of Theorem \ref{main-th4}.

\lemma{le7-1}
{For any $r\in \brk{n+1}$ and partition $\{U_1,U_2\}$ of
$\brk{n+2}\setminus \{r,r+1\}$,
\begin{equation}\label{eq7-1}
\tutte^c_{CT(U_1,\{r\},U_2)}(y)=
\tutte^c_{CT(U_1\cup\{r\},U_2\cup \{r+1\})}(y)
-
\tutte^c_{CT(U_1\cup\{r+1\},U_2\cup \{r\})}(y).
\end{equation}
}

\proof Let $G$ denote the graph
$CT(U_1\cup\{r\},U_2\cup \{r+1\})$ and $e$ the edge
in $G$ joining $r$ and $r+1$.
As $G$ is a simple graph and $e$ is not a loop,
by (\ref{eq2}) and the definition of $\tutte^c_G(y)$,
\begin{equation}\label{eq7-2}
\tutte^c_{G}(y)=\tutte^c_{G\backslash e}(y)+\tutte^c_{G\slash e}(y).
\end{equation}
Observe that $G\backslash e$ and $G\slash e$ are isomorphic
to $CT(U_1\cup\{r+1\},U_2\cup \{r\})$
and $CT(U_1,\{r\},U_2)$ respectively.
Thus, the result follows.
\proofend

\lemma{le7-2}
{Let $G$ be a
 complete tiered graph with two tiers
and $V(G)=\brk{n+2}$.
For $i=1,2$,
let $a_i(G)$ be the cardinality of the set
$\{r\in \brk{n+1}: r\in V_i(G), r+1\in V_{3-i}(G)\}$.
If $G$ is connected, then $a_1(G)-a_2(G)=1$.
}

\proof
Let $\sets_G$ denote the unique  partition
$\{S_1,S_2,\cdots,S_k\}$ of $V(G)$ determined by the following
conditions:
\begin{itemize}
\item[(a).] for each $i\in \brk{k}$, $S_i\subseteq V_j(G)$ holds
for some $j\in \brk{2}$; and
\item[(b).] each $S_i$ is a maximal set of consecutive integers.
\end{itemize}
For example, if $V_1(G)=\{1,3,4,6,7\}$ and
$V_2(G)=\{2,5,8,9,10,11\}$, then $\sets_G=\{S_i:1\le i\le 6\}$,
where $S_1=\{1\}, S_2=\{2\}$, $S_3=\{3,4\}$, $S_4=\{5\}$,
$S_5=\{6,7\}$ and $S_6=\{8,9,10,11\}$.

 Assume that $\min(S_1)<\min(S_2)
<\cdots<\min(S_k)$,
where $\min(S_t)$ is the minimum
number in $S_t$.  It is easily seen that
\begin{enumerate}
\item[(a').] if $S_1\subseteq V_2(G)$, then all vertices in $S_1$ are isolated in $G$;
and

\item[(b').] if $S_1\subseteq V_1(G)$ and $k$ is odd, then $S_k\subseteq V_1(G)$,
implying that all vertices in $S_k$ are isolated in $G$.
\end{enumerate}
As $G$ is connected,  we  have $S_1\subseteq V_1(G)$ and $k=2s$ for some positive integer $s$. This implies that $a_1(G)=s$ and $a_2(G)=s-1$, completing the proof.
\proofend

Now we are ready to prove Theorem~\ref{main-th4}.

\vspace{0.2 cm}

\noindent {\it Proof of Theorem~\ref{main-th4}}:
Let $n=p_1+p_2$.
By Theorem~\ref{main-th}, $P_{(1,p_1,p_2)}(q)=P_{(p_1,1,p_2)}(q)$.
By (\ref{eq1-4}),
\begin{equation}\label{eq7-3}
P_{(p_1,1,p_2)}(q)
=\sum_{G\in \CT_{(p_1,1,p_2)}}\tutte^c_G(q)
=\sum_{r=1}^{n+1}\sum_{G\in \CT_{(p_1,1,p_2)}\atop V_2(G)=\{r\}}
\tutte^c_{G}(q).
\end{equation}
For $1\le r\le n+1$, let
$\setp_{r}$ denote the set of  partitions
$\{U_1,U_2\}$ of $\brk{n+2}\setminus \{r,r+1\}$ with $|U_i|=p_i$
for $i=1,2$.
By (\ref{eq7-3}) and Lemma~\ref{le7-1},
\begin{eqnarray}\label{eq7-4}
P_{(p_1,1,p_2)}(q)
&=&\sum_{r=1}^{n+1}
\sum_{\{U_1,U_2\}\in \setp_{r}}
\left (\tutte^c_{CT(U_1\cup\{r\},U_2\cup \{r+1\})}(q)
-
\tutte^c_{CT(U_1\cup\{r+1\},U_2\cup \{r\})}(q)\right)
\nonumber \\
&=&\sum_{G\in \CT^c_{(p_1+1,p_2+1)}} (a_1(G)-a_2(G)) \tutte^c_{G}(q).
\end{eqnarray}
By  Lemma~\ref{le7-2},  we have   $a_1(G)-a_2(G)=1$.
Hence, Theorem~\ref{main-th4} follows from (\ref{eq7-4})
and (\ref{eq1-4}).
\proofend

\label{r1-5}

\section*{Acknowledgments}
The authors would like to thank the reviewers for their helpful comments and suggestions.
The second author was supported by
the National Natural
Science Foundation of China (11671366 and 12071440).

\end{document}